\documentclass{gtmon_a}
\pdfoutput=1

\usepackage[all]{xy}

\proceedingstitle{The Zieschang Gedenkschrift}
\conferencestart{5 September 2007}
\conferenceend{8 September 2007}
\conferencename{Conference in honour of Heiner Zieschang}
\conferencelocation{Toulouse, France}

\editor{Michel Boileau}
\givenname{Michel}
\surname{Boileau}

\editor{Martin Scharlemann}
\givenname{Martin}
\surname{Scharlemann}

\editor{Richard Weidmann}
\givenname{Richard}
\surname{Weidmann}

\title{On invertible generating pairs of fundamental groups of graph
manifolds}

\author{Michel Boileau}
\givenname{Michel}
\surname{Boileau}
\address{Laboratoire \'Emile Picard\\
CNRS UMR 5580\\
Universit\'e Paul Sabatier\\\newline
118 Route de Narbonne\\
F-31062 Toulouse Cedex 4\\
France\vspace{3pt}\\\newline
Department of Mathematics and the Maxwell Institute of
Mathematical Sciences\\\newline
Heriot-Watt University\\
Riccarton\\
Edinburgh\\
EH14 4AS\\
Scotland}
\email{boileau@picard.ups-tlse.fr}
\urladdr{}

\author{Richard Weidmann}
\givenname{Richard}
\surname{Weidmann}
\email{R.Weidmann@ma.hw.ac.uk}
\urladdr{}

\dedicatory{Dedicated to Heiner Zieschang, our teacher, mentor and friend}

\volumenumber{14}
\issuenumber{}
\publicationyear{2008}
\papernumber{7}
\startpage{105}
\endpage{128}

\doi{}
\MR{}
\Zbl{}

\arxivreference{}


\received{21 October 2006}
\revised{}
\accepted{5 February 2007}
\published{29 April 2008}
\publishedonline{29 April 2008}
\proposed{}
\seconded{}
\corresponding{}
\version{}

\makeatletter
\def\cnewtheorem#1[#2]#3{\newtheorem{#1}{#3}
\expandafter\let\csname c@#1\endcsname\c@Theorem}
\makeatother

\let\xysavmatrix\xymatrix
\def\xymatrix{\disablesubscriptcorrection\xysavmatrix}
\AtBeginDocument{\let\bar\wbar\let\tilde\wtilde\let\hat\what}

\makeop{rank}
\makeop{Stab}
\makeop{Aut}
\makeop{Fix}
\newcommand{\PSL}{\mathit{PSL}}

\newtheorem{Theorem}{Theorem}
\cnewtheorem{Lemma}[Theorem]{Lemma}
\cnewtheorem{Corollary}[Theorem]{Corollary}
\cnewtheorem{Proposition}[Theorem]{Proposition}
\theoremstyle{definition}
\cnewtheorem{Remark}[Theorem]{Remark}
\cnewtheorem{Question}[Theorem]{Question}
\makeautorefname{Theorem}{Theorem}
\makeautorefname{Lemma}{Lemma}
\makeautorefname{Corollary}{Corollary}
\makeautorefname{Proposition}{Proposition}
\makeautorefname{Question}{Question}

\begin{document}

\begin{htmlabstract}
We study invertible generating pairs of fundamental groups of graph
manifolds, that is, pairs of elements (g,h) for which the map g &rarr;
g<sup>-1</sup>, h &rarr; h<sup>-1</sup> extends to an automorphism. We show in
particular that a graph manifold is of Heegaard genus 2 if and only if
its fundamental group has an invertible generating pair.
\end{htmlabstract}

\begin{abstract}
We study invertible generating pairs of fundamental groups of graph
manifolds, that is, pairs of elements $(g,h)$ for which the map $g\mapsto
g^{-1}$, $h\mapsto h^{-1}$ extends to an automorphism. We show in
particular that a graph manifold is of Heegaard genus 2 if and only if
its fundamental group has an invertible generating pair.
\end{abstract}
\begin{asciiabstract}
We study invertible generating pairs of fundamental groups of graph
manifolds, that is, pairs of elements (g,h) for which the map g -->
g^{-1}, h --> h^{-1} extends to an automorphism. We show in
particular that a graph manifold is of Heegaard genus 2 if and only if
its fundamental group has an invertible generating pair.
\end{asciiabstract}

\maketitle

\section{Introduction}

Let $M$ be a closed orientable irreducible 3--manifold. We say that
$M$ is a 2--generated manifold if $G=\pi_1(M)$ can be generated by two
elements. Thurston asked the following question (see also Kirby 
\cite[Problem~3.15A]{Ki}):

\begin{Question}\label{q1}
Given a generating pair $(g,h)$ for $G = \pi_1(M)$, is there an 
automorphism $\alpha\co G\to G$ such that $\alpha(g)=g^{-1}$ and $\alpha(h)=h^{-1}$?
\end{Question}

We call such an automorphism $\alpha$ an {\em inversion} of $\pi_1(M)$,
and if it exists we say that the generating pair $(g,h)$ is {\em
invertible}.  Note that  3--manifolds of Heegaard genus $2$ always have inversions. We say that $\alpha$ is a {\em
hyperelliptic inversion} if $\alpha$ is induced by the hyperelliptic involution of a  
genus $2$  Heegaard splitting. 

\begin{Question}\label{q2}
Is any inversion of $\pi_1(M)$ hyperelliptic?
\end{Question}

For 2--generated  hyperbolic manifolds \fullref{q1} has a positive answer by Jorgensen's 
construction (see Adams \cite{Ad} and Thurston \cite[Chapter~5]{Thu})
while \fullref{q2} 
is still open. In this article we consider the case of graph manifolds. Our main result is the following:

\begin{Theorem}\label{main}  A closed orientable graph manifold $M$ has Heegaard genus two if and only if $\pi_1(M)$ admits an inversion.
\end{Theorem}

As there are graph manifolds and even Seifert manifolds with $2$--generated
fundamental group that do not admit a genus $2$ Heegaard splitting (see
Boileau and Zieschang \cite{BZ} and Weidmann \cite{We}) we get  a negative
answer to  \fullref{q1}:

\begin{Corollary}
There exist graph manifolds with 2--generated fundamental group whose
fundamental groups have no invertible generating pair.
\end{Corollary}

In the case of Seifert manifolds we have a stronger result, namely: 

\begin{Theorem}\label{thm:seifert}
Let $M$ be a closed orientable Seifert $3$--manifold. Then any inversion
of $\pi_1(M)$ is  hyperelliptic.
\end{Theorem}

For Seifert $3$--manifolds this gives a positive answer to \fullref{q2}. We
believe that, using the same methods, this result can be extended more
generally to graph manifolds.

Even in the case of Seifert manifolds, it turns out that non-invertible
generating pairs are rather plentiful and also  many 3--manifolds of genus
$2$ have such generating pairs. For example Heegaard  genus 2 Seifert
manifolds with base $S^2$ and three singular fibers of order $\{2, 3, r
\}$ such that $\gcd(6,r)=1$ and $r \geq 7$ have infinitely many generating
pairs which are not invertible. The fundamental group of some of these
examples  are 2--generated subgroups of $\widetilde{\PSL}(2,\mathbb{R})$.

It should be remarked also that \fullref{thm:seifert} does not imply
that any invertible generating pair of the fundamental group of a Seifert
manifold is geometric, where we say that a generating pair is geometric
if it corresponds to the spine of a handlebody of a genus $2$ splitting
of the manifold. In fact there are examples of invertible generating
pairs that are seemingly unrelated to the geometric generating pair
that induces the same inversion. The situation seems most mysterious
for the Brieskorn homology sphere wich is a Seifert manifold with base
space $S^2(2,3,7)$. It turns out that its fundamental group admits
an invertible generating pair which is not geometric after Nielsen
equivalence and passing to roots.

\section{Preliminaries}

Let $G$ be a group and $(g,h)$ be a generating pair of $G$. If the map $\bar\alpha(g)=g^{-1}$ and $\bar \alpha(h)=h^{-1}$ extends to an automorphism $\alpha$ of the group $G$ then we say that the generating pair
$(g,h)$ induces an {\em inversion} $\alpha$ of $G$, sometimes we call $\alpha$ the inversion {\emph relative $(g,h)$}. We then further say that the generating pair $(g,h)$ is {\em invertible.}

It is clear that for a given group some generating pair might induce an inversion while another does not. In the following however  we will see that this property is an invariant of the Nielsen equivalence class.

Recall that two tuples $T=(g_1,\ldots ,g_n)$ and $T'=(h_1,\ldots ,h_n)$ are {\em elementary equivalent} if  one of the following holds:

\begin{enumerate}
\item $h_i=g_{\sigma(i)}$ for $1\le i\le n$ and some $\sigma\in S_n$.
\item $h_1=g_1^{-1}$  and $h_i=g_i$ for $i\ge 2$.
\item $h_1=g_1g_2$ and $h_i=g_i$ for $i\ge 2$.
\end{enumerate}

We further say that two tuples $T$ and $T'$ are {\em Nielsen equivalent} if there is a sequence of tuples $T=T_0,T_1,\ldots ,T_{k-1},T_k=T'$ of tuples such that $T_{i-1}$ and $T_i$ are elementary equivalent for $1\le i\le k$.

\begin{Lemma}\label{lem:nielsen} Let $G$ be a group and $(g,h)$ and $(g',h')$ be Nielsen-equivalent generating pairs. Then $(g,h)$ is invertible if and only if $(g',h')$ is invertible. 
\end{Lemma}

\begin{proof}
It clearly suffices to verify the assertion for pairs $(g,h)$ and $(g',h')$ that are elementary equivalent. If they are related as in (1) or (2) above  then the assertion is obvious as the maps $\bar \alpha$ are defined the same way.

Thus we can assume that $g'=gh$ and $h'=h$. Suppose that $\alpha_1\co G\to
G$ is an inversion that extends the map $\bar\alpha_1\co \{g,h\}\to G$
given by $\bar\alpha_1(g)=g^{-1}$ and $\bar\alpha_1(h)=h^{-1}$. We have to
show that the map $\bar\alpha_2\co \{gh,h\}\to G$ given by
$\bar\alpha_2(gh)=h^{-1}g^{-1}$ and $\bar\alpha_2(h)=h^{-1}$ extends to an
automorphism $\alpha_2\co G\to G$, it is easily verified that this map $\alpha_2$ can be defined as $\alpha_2(x)=h^{-1}\alpha_1(x)h$ for all $x\in G$. This map is clearly an automorphism of $G$ and its restriction to $\{gh,h\}$ is $\bar\alpha_2$.
\end{proof}

Occasionally we will change a generating pair in such a way that we replace one generator with a root. This operation also turns out to be well behaved:

\begin{Lemma}\label{lem:root} Let $G$ be a group, $g,h\in G$ and $n\ge 2$ such that $(g,h)$ and $(g^n,h)$ are generating pairs. Suppose that $(g,h)$ induces an inversion. Then $(g^n,h)$ induces the same inversion.
\end{Lemma}

\begin{proof}
This is immediate as $(g^n,h)$ is a generating pair
and the induced inversion $\alpha$ must map $g^n$ to
$\alpha(g^n)=\alpha(g)^n=\left(g^{-1}\right)^n=g^{-n}$.
\end{proof}

The following corollary is  straightforward consequence of the proofs
of \fullref{lem:nielsen} and \fullref{lem:root}.

\begin{Corollary}\label{cor:extnielsen}
Let $(g',h')$ be an invertible generating pair of a group $G$ obtained
from a generating pair $(g,h)$ by Nielsen equivalence and replacing $g$
and/or $h$ with roots. Then $(g,h)$ is also invertible and the inversions
$\alpha$ induced by $(g,h)$ and  $\alpha'$ induced by $(g',h')$ give
the same element in
$Out(G)$ (the outer automorphism group of $G$).
\end{Corollary}

 If $G$ is the fundamental group of a 3--manifold $M$ then there are
situations when we can guarantee that a generating pair induces an
inversion. If  $M$ is a hyperbolic 3--manifold then the involution exists
for any generating pair, this result is due to Jorgensen, see Adams \cite{Ad} or
Thurston \cite[Chapter~5]{Thu}.

Another situation where we can guarantee the existence of the inversion
is the case of geometric generating pairs: Suppose that the manifold $M$
is of Heegaard genus $2$, that is, that $M=V_1\cup_SV_2$ where $V_1$ and $V_2$
are handlebodies of genus $2$ glued together along their boundaries. Now
$\pi_1(V_i)$ is free of rank $2$ and the homomorphism $i_*\co \pi_1(V_i)\to
\pi_1(M)$ induced by the inclusion map is surjective. This implies that
any generating pair $(g,h)$ of $\pi_1(V_i)$ gets mapped onto a generating
pair  $(i_*(g),i_*(h))$ of $\pi_1(M)$. We call generating pairs that
occur this way {\em geometric generating pairs}.

It has been shown by Birman and Hilden \cite{BH} that any
geometric generating set induces an inversion. To see this it suffices
to show that the hyperelliptic involution of the handlebody of genus
$2$ induces an automorphism of the fundamental group that maps the
generators to their inverses and that this involution can be extended to
the complementary handlebody. We call an inversion of $\pi_1(M)$ that is
obtained as above a {\em hyperelliptic inversion} of $\pi_1(M)$. As any
generating pair of $\pi_1(M)$ that is Nielsen equivalent to a generating
pair that corresponds to the spine of handlebody in $M$ also corresponds
to a spine of the same handlebody we have the following observation.

\begin{Lemma}\label{lem:geom}  Let $G$ be the fundamental group of a $3$--manifold $M$ and $(g,h)$ and $(g',h')$ be Nielsen-equivalent generating pairs. Then $(g,h)$ is geometric if and only if $(g',h')$ is geometric. \hfill$\Box$
\end{Lemma}

\section{Compact 2--orbifolds with $2$--generated fundamental group}  

In this section we study generating pairs of Fuchsian groups, that is, of
fundamental groups of $2$--orbifolds. Throughout this article we assume
that all 2--orbifolds have no reflection boundary lines. As all compact
bad 2--orbifolds have finite fundamental groups it follows that any
compact 2--orbifold is either Euclidean or hyperbolic or has finite
fundamental group. We need to understand generating pairs of Euclidean and
hyperbolic 2--orbifolds.

\begin{Lemma}\label{euclidean} Let $\mathcal O$ be a compact euclidean
2--orbifold with $2$--generated fundamental group. Then one of the following holds:
\begin{enumerate}
\item ${\mathcal O}=T^2$, $G=\pi_1({\mathcal O})=
  \langle a,b\,|\,[a,b]\rangle$ and every pair of generators of $G$
  is Nielsen equivalent to $(a,b)$.
\item ${\mathcal O}=$KB, $G=\pi_1({\mathcal O})=
  \langle a,b\,|\,aba^{-1}b\rangle$ and every pair of generators of $G$
  is Nielsen equivalent to $(a,b)$.
\item ${\mathcal O}=S^2(p,q,r)$,
  $1=\frac{1}{p}+\frac{1}{q}+\frac{1}{r}$ and $G=\pi_1(\mathcal{O})=
  \langle s_1,s_2,s_3\,|\,s_1^p,s_2^q,s_3^r,s_1s_2s_3\rangle\allowbreak=
  \langle s_1,s_2\,|\,s_1^p,s_2^q,(s_1s_2)^r\rangle$. Moreover any
  generating pair $(g,h)$ of $G$ is Nielsen-equivalent to $(s_1,s_2)$.
\item  ${\mathcal O}=P^2(2,2)$, that is, $G=\pi_1({\mathcal O})=
  \langle x,s_1,s_2\,|\,s_1^2,s_2^2,x^2s_1s_2\rangle$. Moreover any
  generating pair is of $G$ is Nielsen-equivalent to $(x,s_1)$.
\end{enumerate}
\end{Lemma}

\begin{proof}
Note first that every compact Euclidean 2--orbifold
with 2--generated fundamental group occurs in this list as the only
remaining compact Euclidean 2--orbifold  $S^2(2,2,2,2)$ has a fundamental group of rank $3$. In the first two cases the assertion is obvious.

Suppose next  that ${\mathcal O}=S^2(p,q,r)$ with $1=\frac1p+\frac1q+\frac1r$. This clearly implies that  $(p,q,r)$ is $(3,3,3)$, $(2,3,6)$ or $(2,4,4)$. Note first that any rotation $r$ that is a member of a generating pair of $G$ must generate a maximal cyclic subgroup as otherwise $G/\langle\langle r\rangle\rangle$ is non-cyclic. As any maximal finite cyclic subgroup is of order $2$, $3$, $4$ or $6$ this implies that any such rotation must be conjugate to $s_i^{\pm 1}$ for some $i$. Let now $(g,h)$ be a generating pair. We can assume that $g$ is a rotation as two translations necessarily generate a subgroup of type $\mathbb Z^2$. Thus we can assume that $(g,h)=(s_i,h)$ for some $i$. Possibly after replacing $h$ with $hg$, we can assume that $h$ is also a rotation as otherwise $\langle g,h\rangle=\langle h,hg\rangle$ is free Abelian. It follows that  $h$ is  conjugate to $s_j$ for some $j$, again possible after replacing $h$ with $h^{-1}$. 

Note first that if $g$ and $h$ are both of order $6$, that is, if $(p,q,r)=(2,3,6)$ and $g$ and $h$ are conjugate to $s_3$, then we can replace $g$ by $gh$ and obtain a rotation of order $3$.
It is now obvious that $g$ and $h$ are the standard generating set of a Euclidean triangle group $\pi_1(S^2(l,m,n))$ where $l$ is the order of $g$, $m$ is the order of $h$ and $n$ is the order of $gh$. If $g$ and $h$ generate $G$ then we must have $(l,m,n)$ is a permutation of $(p,q,r)$. The assertion follows easily.

We conclude by looking at the case ${\mathcal O}=P^2(2,2)$. Note that
in this case $G$ splits as an amalgamated product over the cyclic
subgroup generated by $x^2$, namely $G=\langle x\rangle*_{\langle
x^2=s_2s_1\rangle}\langle s_1,s_2\rangle$. Now any generating
pair of a $2$--generated amalgamated product is Nielsen equivalent
to a generating pair $(g',h')$ such that a power of $g'$ lies in
the amalgam. Note first that $g'$ cannot lie in the amalgam as
$G/\langle \langle x^2\rangle\rangle$ is an infinite dihedral group
and therefore not cyclic. As $\langle s_2s_1\rangle$ is isolated in
$\langle s_1,s_2\rangle$ this implies that (possibly after conjugation)
$g'\in\langle x\rangle-\langle x^2\rangle$, that is, that $g'=x^{2l+1}$
for some $l\ge 0$. Suppose first that $l\ge 1$. After passing to
the quotient $$G':=G/\langle\langle {g'}^2\rangle\rangle=\langle
x|x^{4l+2}\rangle*_{\langle x^2=s_2s_1\rangle}\langle
s_1,s_2|s_1^2,s_2^2,(s_1s_2)^{2l+1}\rangle$$ and replacing $g'$ and
$h'$ with their images in $G'$ we see that $g'$ does not lie in any
conjugate of the amalgam of the new amalgamated product. It therefore
follows from Kapovich and Weidmann
\cite[Lemmas~2.1 and~2.2]{KW} that after replacing
$h'$ with $h''={g'}^{z_1}h'{g'}^{z_2}$ for some $z_1,z_2\in\mathbb Z$
a power of $h''$ also lies in $\langle x\rangle$.
This however cannot be as $g'$ and $h''$ both lie in the kernel of the
quotient map $G'\to G'/\langle\langle x\rangle\rangle\cong \mathbb Z_2$,
thus the original elements $g$ and $h$ could not have generated $G$. Thus
we can assume that $g'=x$. A simple normal form argument then shows
that after Nielsen equivalence either $h'\in\langle s_1,s_2\rangle$
or that $(g',h')$ do not generate $G$. The assertion now follows
easily.
\end{proof}

We lastly deal with hyperbolic 2--orbifolds. In the case of 2--orbifolds
with orientable base space this situation has been completely discussed
by Fine and Rosenberger~\cite{FR}, the remaining cases turn out be
relatively easy.

\begin{Lemma}\label{hyporb} Let $\mathcal O$ be a compact hyperbolic
2--orbifold such that $G=\pi_1({\mathcal O})$ is $2$--generated. Then one of the following holds:
\begin{enumerate}
\item ${\mathcal O}=T^2(p)$, $G=\langle a,b\,|\,[a,b]^p\rangle$ and every pair of generators of $G$ is Nielsen equivalent to $(a,b)$.
\item ${\mathcal O}=$KB(p), $G=\langle a,b\,|\,(aba^{-1}b)^p\rangle$ and every pair of generators of $G$ is Nielsen equivalent to $(a,b)$.
\item  ${\mathcal O}=P^2(p,q)$ with $p\ge 3$ and $q\ge 2$. Thus $G=\langle x,s_1,s_2\,|\,s_1^p,s_2^q,x^2s_1s_2\rangle$. Moreover any generating pair is Nielsen equivalent to $(x,s_1)$ unless  $p=2$ and $q$ is odd in which case any generating set is Nielsen equivalent to either $(x,s_1)$ or $(x^2,xs_1)$.
\item ${\mathcal O}=S^2(2,2,2,2l+1)$, $G=\langle s_1,\ldots ,s_4\,|\, s_1^2,s_2^2,s_3^2,s_4^{2l+1},s_1s_2s_3s_4\rangle$ and any generating pair is Nielsen equivalent to $(s_1s_2,s_1s_3)$.
\item ${\mathcal O}=S^2(p,q,r)$ with  $1>\frac{1}{p}+\frac{1}{q}+\frac{1}{r}$ and $G=\langle s_1,s_2\,|\,s_1^p,s_2^q,(s_1s_2)^r\rangle=\langle s_1,s_2,s_3\,|\,s_1^p,s_2^q,s_3^r,s_1s_2s_3\rangle$. Moreover any generating pair $(g,h)$ is either Nielsen equivalent to a pair of type $(s_i^{n_1},s_j^{n_2})$ with $1\le i,j\le 3$ and $n_1,n_2\in\mathbb N$ or one of the following holds:
\begin{enumerate}
\item $(p,q,r)=(2,3,r)$ with $\gcd(6,r)=1$ and $(g,h)$ is Nielsen equivalent to $(g',h')=(s_1s_2s_1s_2^{-1},s_2^{-1}s_1s_2s_1)$.
\item $(p,q,r)=(2,4,r)$ with odd $r$ and $(g,h)$ is Nielsen equivalent to $(g',h')=(s_1s_2^2,s_2^{-1}s_1s_2^{-1})$.
\item $(p,q,r)=3,3,r$ with $\gcd(3,r)=1$ and $(g,h)$ is Nielsen equvalent to $(g',h')=(s_1s_2^{-1},s_2^{-1}s_1)$.
\item $(p,q,r)=(2,3,7)$ and $(g,h)$ is Nielsen equivalent to  $(g',h')=$\hfill\break $(s_1s_2^{-1}s_1s_2s_1s_2^{-1}s_1s_2^{-1}s_1s_2,s_2^{-1}s_1s_2s_1s_2^{-1}s_1s_2s_1s_2s_1)$.
\end{enumerate}
\end{enumerate}
\end{Lemma}

\begin{proof}
Note first that the list contains all compact
hyperbolic $2$--orbifolds with $2$--generated fundamental group by
Zieschang~\cite{Zie1} and Peczynski, Rosenberger and Zieschang~\cite{PRZ}.
The assertion of the lemma then follows from Fine and Rosenberger~\cite{FR} except in situation (2) and (3).

If ${\mathcal O}=KB(p)$ then $G=\langle a,b,s\,|\,
s^p,aba^{-1}=sb^{-1}\rangle$, that is, $G$ is a HNN--extension of the group
$\langle b,s\rangle\cong\mathbb Z*\mathbb Z_p$ with stable letter $a$ and
associated subgroups $\langle b\rangle$ and $\langle sb^{-1}\rangle$.
After Nielsen equivalence and conjugation we can assume that one
generator, say $g$, lies in $\langle b\rangle$ (see, for example, Kapovich and
Weidmann~\cite{KW}) and it is not difficult to verify that $g$ cannot be a
proper power of $b$ as  $G/\langle \langle b^n\rangle\rangle$ is not
cyclic for $n\ge 2$. Thus we can assume that $g=b$. It is further a simple
calculation in normal forms to see that after replacing $h$ with
$g^{z_1}hg^{z_2}$ for appropriate $z_1,z_2\in\mathbb Z$ we can assume that
$h=a$ which proves the assertion, the essential step is to apply
\cite[Lemma~5]{BW}.

Let ${\mathcal O}=P^2(p,q)$. As in the proof of \fullref{euclidean} we see
that we can assume that either $g=x$ or that $g=x^2$. If $g=x$ then we
argue as in the proof of \fullref{euclidean} and get $h=s_1$. If $g=x^2$
then it is easy to see that $h$ must be of type $x\bar h$
\cite[Lemma~2]{BW} and that $h^{-1}gh=\bar h^{-1}x^{-1}x^2x\bar h=\bar
h^{-1}x^2\bar h$ and $x^2$ must generate $\langle s_1,s_2\rangle$. The
last assertion holds as if $\bar h^{-1}x^2\bar h$ and $x^2$ do not
generated $\langle s_1,s_2\rangle$ then it is easily verified that they
generate a free group in which case the original generators cannot have
generated $G$ \cite{BW}. The assertion then follows from the fact that any
orbifold of type $D(p,q)$ whose fundamental group is generated by two
conjugates of the boundary curve is of type $D(2,2l+1)$ (see Rost and
Zieschang~\cite{RZ}).
\end{proof}

\begin{Remark} In all cases of (5) in the above lemma we have that
$[g',h']$ is a power of $s_1s_2$, namely $[g',h']=(s_1s_2)^6$ in case (a),
$[g',h']=(s_1s_2)^4$ in case  (b), $[g',h']=(s_1s_2)^3$ in (c) and
$[g',h']=(s_1s_2)^4$ in case (d). However the cases (a)--(c) are essentially different from  case (d) in that in the first three cases the identity $[g',h']=(s_1s_2)^r$ holds in $\langle s_1,s_2|s_1^p,s_2^q\rangle$ while in case (d) the relation $(s_1s_2)^r$ is needed.\end{Remark}

We conclude this section by studying the group $\pi_1(S^2(2,3,7))$ in somewhat more detail.

\begin{Lemma}\label{inner237} Let $G=\pi_1(S^2(2,3,7))=\langle s_1,s_2\,|\,s_1^2,s_2^3,(s_1s_2)^7\rangle$. Then the pair $(g,h)$ given in 5 (d) above induces an inversion $\alpha$ of $G$. 
Moreover $\alpha$ is the inner automorphism given by $$w\mapsto (s_2s_1s_2^{-1})w(s_2s_1s_2^{-1})^{-1}.$$\end{Lemma}

\begin{proof}
It clearly suffices to show that
$$(s_2s_1s_2^{-1})g(s_2s_1s_2^{-1})^{-1}=g^{-1} \quad\text{and}\quad
(s_2s_1s_2^{-1})h(s_2s_1s_2^{-1})^{-1}=h^{-1}.$$
This is a simple computation.
\end{proof}

\section{Seifert manifolds}\label{sec:seifert}

In this section we establish \fullref{thm:seifert} for Seifert
manifolds with 2--generated fundamental group except for those that
fibre over the base space $P^2(p,q)$ or $S^2(2,2,2,2l+1)$. Note first
that any orientable Seifert manifold with spherical base space is
either a lens space or can be fibered over a base space of type
$S^2(p,q,r)$ with $\frac{1}{p}+\frac{1}{q}+\frac{1}{r}>1$. In the case
of a lens space there is nothing to show and in the latter case any
generating pair projects onto a pair of rotations in
$\pi_1(S^2(p,q,r))\le SO(3)$. Together with the previous section this
implies that for any orientable Seifert manifold $M$ with generating
pair $(g,h)$ and base space different from $S^2(2,2,2,2l+1)$ and
$P^2(p,q)$ one of the following holds:

\begin{enumerate}
\item After Nielsen equivalence we can assume that $g$ and $h$ map to elliptic elements of the base group.
\item $M$ is Seifert fibered with base of type $T$, $T(p)$, $KB$ or $KB(p)$.
\item $M$ is fibered over a hyperbolic base space of type $S^2(p,q,r)$ with $(p,q,r)$ as in (5) of \fullref{hyporb} and $(g,h)$ projects on the respective generating set of $\pi_1(S^2(p,q,r))$.
\end{enumerate}

We now deal with each of these three cases.

\begin{Lemma} Let $M$ be a Seifert manifold with $2$--generated fundamental group. Any generating pair $(g,h)$ which projects to elliptic elements  in the base group can be inverted by a hyperelliptic inversion of  
$\pi_1(M)$.
\end{Lemma}

\begin{proof}
The assumption that the images $\pi(g)$ and $\pi(h)$ of the two generators are elliptic implies that $M$ is Seifert fibered over a base of type $S^2(p,q,r)$.  After possibly replacing $g$ and $h$ with roots we may assume that the element  $f$ corresponding to the fibre is a power of $g$ and $h$. Then $g$ and $h$ are inverted by the inversion induced by the Montesinos involution
$\tau$ on $M$ (see Montesinos~\cite{Mon1,Mon2}). This involution is geometric, reverses the orientation of the base $S^2(p,q,r)$ while fixing each singular point and reverses the orientation of the generic fiber $f$. The quotient map $M \to M/\tau$ 
 is a $2$--fold covering of $S^3$ branched over a Montesinos link $L$ with
$3$ bridges by Boileau and Zieschang~\cite{BZ2}. Therefore $\tau$ is
induced by the hyperelliptic involution of a genus $2$ Heegaard splitting
of $M$ obtained by lifting the $3$--bridge presentation of $L$.
\end{proof}

\begin{Lemma} Let $M$ be a Seifert manifold with $2$--generated fundamental group that fibers over $T^2$, $T^2(p)$, $KB$ or $KB(p)$. Then any generating pair of  $\pi_1(M)$ is geometric. In particular it can be inverted by a hyperelliptic inversion.
\end{Lemma}

\begin{proof}
We consider $T^2$ as $T^2(1)$ and $KB$ as $KB(1)$.
Recall that the fundamental group is either of the form
$\langle a,b,f\,|\,[a,f],b],f],[a,b]^p=f^e\rangle$ or $\langle
a,b,f\,|\,afa^{-1}=f^{-1},bfb^{-1}=f^{-1},(a^2b^2)^p=f^e\rangle$ with
$e=1$ as otherwise $\pi_1(M)$ is not $2$--generated (see Boileau and
Zieschang~\cite{BZ}).

In both cases it is easy to see that $(a,b)$ is a geometric generating set
corresponding to a horizontal Heegaard splitting. It is futher clear that
any pair $(af^l,bf^k)$ is also a generating set as $[af^l,bf^k]=[a,b]$
and $(af^l)^2(bf^l)^2=a^2b^2$. The map $a\mapsto af^l$ and $b\mapsto bf^k$
clearly extends to an isomorphism which is induced by homeomorphism. Thus
$(gf^l,bf^k)$ is also geometric.  \end{proof}

We will (sometimes implicitly) use the following simple observation.

\begin{Lemma}\label{intersectionwithf} Let $M$ be a Seifert manifold
fibered over $S^2(p,q,r)$ and $$\pi\co \pi_1(M)\to \langle s_1,s_2\,|\,s_1^p,s_2^q,(s_1s_2)^r\rangle$$ be the projection. Let $U=\langle g_1,g_2,x_1,\ldots ,x_k\rangle\le \pi_1(M)$ such that $\pi(g_i)=s_i$ for $i=1,2$ and $x_i\in\langle f\rangle=\ker \pi$ for $1\le i\le k$. Then $U\cap\langle f\rangle=\langle x_1,\ldots ,x_k,g_1^p,g_2^q,(g_1g_2)^r\rangle$.
\end{Lemma}

\begin{Proposition}\label{case23r} Let $M$ be a Seifert manifold with hyperbolic base space $S^2(2,3,r)$ where $\gcd(6,r)=1$, $S^2(2,4,r)$ with odd $r$ or $S^2(3,3,r)$ with $\gcd(r,3)=1$. 

Let $(g,h)$ be a generating pair of $\pi_1(M)$ that projects onto the
generating pair of the base group spelt out in \fullref{hyporb} (5) (a)--(c). Then the following hold:
\begin{enumerate}
\item If $(g,h)$ is invertible then $(g,h)$ is geometric and corresponds to a horizontal Heegaard splitting of genus $2$.
\item There exist non-invertible generating pairs of this type.
\end{enumerate}
\end{Proposition}

\begin{proof}
We give all details in the case of $S^2(2,3,r)$ with $\gcd(6,r)=1$ and then comment on the other cases.

Note that $r=6p+\varepsilon$ with $\varepsilon\in\{-1,1\}$ and $p\ge 1$. Thus we assume that $g=s_1s_2s_1s_2^{-1}f^k$ and $h=s_2^{-1}s_1s_2s_1f^l$ generate the group
\begin{align*}
G&=\langle s_1,s_2,s_3,f\,|\,s_1^2=f,s_2^3=f,s_3^r=f^n,s_1s_2s_3=f^m\rangle \\
&=\langle s_1,s_2,f\,|\,s_1^2=f,s_2^3=f,(s_1s_2)^r=f^{q}\rangle
\end{align*}
where $q=rm-n$.

Suppose now that $(g,h)$ is invertible and let $\alpha$ be the induced inversion. The homeomorphism corresponding to $\alpha$ must be isotopic to the identity as $\alpha$ preserves the fibre. Thus $\alpha$ must be an inner automorphism. It is easily verified that $g$ is conjugate to $g^{-1}$ and $h$ is conjugate to $h^{-1}$ if and only if $k=l=-1$. 

We need some calculations.  Note first that
\begin{multline*}
[g,h]=s_1s_2s_1s_2^{-1}\cdot s_2^{-1}s_1s_2s_1\cdot
s_2s_1^{-1}s_2^{-1}s_1^{-1}\cdot s_1^{-1}s_2^{-1}s_1^{-1}s_2 \\
=s_1s_2s_1s_2^{-2}s_1s_2s_1s_2s_1^{-1}s_2^{-2}s_1^{-1}s_2f^{-1}=(s_1s_2)^6f^{-5}
\end{multline*}
and therefore
\begin{align*}
[g,h]^r=& (s_1s_2)^{6r}f^{-5r}=f^{6q-5r} \\
\text{and}\qquad\qquad c:=& [g,h]^p=(s_1s_2)^{6p}f^{-5p}=(s_1s_2)^{r-\varepsilon}f^{-5p}=(s_1s_2)^{-\varepsilon}f^{q-5p}.
\end{align*}
We define $g':=c^{2\varepsilon}g$, that is, we have
\begin{align*}
g'&=(s_1s_2)^{-2}s_1s_2s_1s_2^{-1}f^{k+2\varepsilon q-10\varepsilon
p} \\
&=s_2^{-2}f^{k+2\varepsilon q-10\varepsilon p}=s_2f^{k-1+2\varepsilon
q-10\varepsilon p} \\
\text{and}\qquad {g'}^3&=f\cdot f^{3k-3+6\varepsilon q-30\varepsilon
p}=f^{3k-2+6\varepsilon q-30\varepsilon p}=f^{3k+3+6\varepsilon
q-5\varepsilon r}.
\end{align*}
We then put $h':={g'}^{-1}c^{-\varepsilon}h$, that is,
\begin{align*}
h'&=s_2^{-1}f^{-k+1-2\varepsilon q+10\varepsilon p}(s_1s_2)f^{-\varepsilon
q+\varepsilon 5p}s_2^{-1}s_1s_2s_1f^l \\
&=s_1f^{-k+l+2-3\varepsilon q +\varepsilon 15p} \\
\text{and}\qquad {h'}^2 &=f^{-2k+2l+5-6\varepsilon
q +\varepsilon 30p}=f^{-2k+2l-6\varepsilon q+5\varepsilon
r}.
\end{align*}
In the case $l=k=-1$ a simple calculation further shows that
$h'g'=c^{-\varepsilon}$ and  it follows that $\langle g,h\rangle=\langle
c,g,h\rangle=\langle g',h'\rangle$. \fullref{intersectionwithf}
then implies that $\langle g,h\rangle\cap \langle f\rangle=\langle
{g'}^3,{h'}^2,(h'g')^r\rangle=\langle f^{6q-5r}\rangle$. As we assume
that $(g,h)$ is a generating set this implies that $6q-5r\in\{-1,1\}$,
in particular $[g,h]$ is the $r^{\rm th}$ root of the fibre, thus there
is a horizontal splitting with respect to which $(g,h)$ is a geometric
generating pair.
It further follows from the above computations that for many other choices of $l$ and $k$ the pair $(g,h)$ is also a generating pair, which must be non-invertible. 

The case of a base space of type $S^2(2,4,r)$ and odd $r$ is completely analogous, again the induced involution must be isotopic to the identity.

In the case of a base space of type $S^2(3,3,r)$ and $r=3p+\varepsilon$
with $\varepsilon\in\{-1,1\}$ the argument is different. We assume
that there is an invertible generating pair $(g,h)$ of
\begin{align*}
G&=\langle s_1,s_2,s_3,f\,|\,s_1^3=f,s_2^3=f^{\eta},
  s_3^r=f^n,s_1s_2s_3=f^m\rangle \\
&=\langle s_1,s_2,f\,|\,s_1^3=f,s_2^3=f^{\eta},(s_1s_2)^r=f^{q}\rangle
\end{align*}
where $\eta\in\{-1,1\}$ and $q=rm-n$ such that
$(g,h)=(s_1s_2^{-1}f^k,s_2^{-1}s_1f^l)$.

As before we see that a power of the commutator $[g,h]$ is a power of $f$.
This implies that  the induced inversion $\alpha$ of $G$ must fix $f$.
Thus  it is induced by a fiber preserving involution of the manifold that
preserves the orientation by Zieschang and Zimmerman~\cite{ZZ}, see also
Zieschang~\cite{Zie2}. This involution cannot be isotopic to the identity as $g=s_1s_2^{-1}f^k$ is not conjugate to $g^{-1}=s_2s_1^{-1}f^{-k}$ for any choice of $k$. This can be seen as $s_1s_2^{-1}$ and $s_2s_1^{-1}$ are not conjugate in the base group. 

It follows that the involution exchanges the two exceptional fibres with invariants $(3,1)$ and $(3,\eta)$, in particular $\eta=1$. Thus the inversion $\alpha$ must map $s_1$ to $s_2$ and $s_2$ to $s_1$ which implies that $k=l=0$. The remainder of the argument is as before.
\end{proof}

For Seifert manifolds $M$ that fibre over the base space $S^2(2,3,7)$, we will now study generating pairs  that project onto the generating pair of the base described  in \fullref{hyporb} (5) (d). Recall that we have 
\begin{align*}
\pi_1(M)&=\langle s_1,s_2,s_3,f\,|\,
s_1^2=f,s_2^3=f,s_3^7=f^n,s_1s_2s_3=f^m\rangle\\
&=\langle s_1,s_2,f\,|\,
s_1^2=f,s_2^3=f,(s_1s_2)^7=f^q\rangle
\end{align*}
with $q:=7m-n$. This shows in
particular that there is a 1--parameter family of Seifert manifolds with base space $S^2(2,3,7)$. The following lemma implies that there are plenty of generating pairs of $\pi_1(M)$ that project onto the generating set of the base group described in \fullref{hyporb} 5(d).

\begin{Lemma}\label{genset237} Let $M$ be an orientable Seifert manifold with base space $S^2(2,3,7)$ and $G=\pi_1(M)=\langle s_1,s_2,f\,|\, s_1^2=f,s_2^3=f,(s_1s_2)^7=f^q\rangle$ with $q$ as above.

Then $G$ is generated by two elements
$$g=s_1s_2^{-1}s_1s_2s_1s_2^{-1}s_1s_2^{-1}s_1s_2f^k
\qquad\text{and}\qquad
h=s_2^{-1}s_1s_2s_1s_2^{-1}s_1s_2s_1s_2s_1f^l$$
iff $\langle f\rangle$
is generated by
$$f^{k+l+5},\quad f^{2k-20q+121},\quad f^{6k+6l+12q-40}\quad\text{and }
f^{-7k-14l-97q+511}.$$
\end{Lemma}
 
\begin{proof}
Note first that
\begin{align*}
gh&=s_1s_2^{-1}s_1s_2s_1s_2^{-1}s_1s_2^{-1}s_1s_2
  f^ks_2^{-1}s_1s_2s_1s_2^{-1}s_1s_2s_1s_2s_1f^l \\
&=s_1s_2^{-1}s_1s_2s_1s_2^{-1}s_1s_2^{-1}s_1^2s_2
  s_1s_2^{-1}s_1s_2s_1s_2s_1f^{k+l} \\
&=s_1s_2^{-1}s_1s_2s_1s_2^{-1}s_1^2s_2^{-1}s_1s_2s_1s_2s_1f^{k+l+1} \\
&=s_1s_2^{-1}s_1s_2s_1s_2^{-2}s_1s_2s_1s_2s_1f^{k+l+2} \\
&=s_1s_2^{-1}(s_1s_2s_1s_2s_1s_2s_1s_2s_1)f^{k+l+1} \\
&=s_1s_2^{-1}s_2^{-1}s_1^{-1}s_2^{-1}s_1^{-1}s_2^{-1}f^{q+k+l+1} \\
&=s_1s_2^{-2}s_1^{-1}s_2^{-1}s_1^{-1}s_2^{-1}f^{q+k+l+1} \\
&=s_1s_2s_1s_2^{-1}s_1^{-1}s_2^{-1}f^{q+k+l-1}.
\end{align*}
Similar arguments show that
$$g^{-1}h^{-1}=s_2s_1s_2^{-1}s_1s_2s_1s_2f^{q-k-l-13}$$
and therefore
\begin{align*}
[g,h]=&ghg^{-1}h^{-1}=(s_1s_2)^4f^{2q-15} \\
\text{and}\qquad\qquad c:=&[g,h]^2=s_1s_2f^{5q-30}.
\end{align*}
Clearly $\langle
g,h\rangle=\langle g,h,c\rangle$ and therefore $\langle
g,h\rangle=\langle g':=c^{-1}gc^{-1},h':=chc,c\rangle$ where
\begin{align*}
g'&=f^{-5q+30}s_2^{-1}s_1^{-1}(s_1s_2^{-1}s_1s_2s_1s_2^{-1}
  s_1s_2^{-1}s_1s_2f^k)f^{-5q+30}s_2^{-1}s_1^{-1} \\
&=s_2s_1s_2s_1s_2^{-1}s_1^{-1}s_2^{-1}f^{-10q+k+60} \\
\text{and}\qquad
h'&=s_1s_2f^{5q-30}(s_2^{-1}s_1s_2s_1s_2^{-1}s_1s_2s_1s_2s_1f^l)
 s_1s_2f^{5q-30} \\
&=s_2s_1s_2^{-1}s_1s_2s_1^{-1}s_2^{-1}f^{10q+l-56}.
\end{align*}
Conjugation with $s_2s_1$ shows that $\langle
g,h\rangle$ is conjugate to $\langle g'',h'',c'\rangle$ with
\begin{align*}
g''&=s_2s_1s_2^{-1}f^{-10q+k+60},\\
h''&=s_2^{-1}s_1s_2f^{10q+l-56} \\
\text{and}\qquad
c'&=s_1^{-1}s_2^{-1}(s_1s_2f^{5q-30})s_2s_1
=s_1s_2^{-1}s_1s_2^{-1}s_1f^{5q-30}.
\end{align*}
We further replace $c'$ with $c''=g''c'h''$ which gives
\begin{align*}
c''&=s_2s_1s_2^{-1}f^{-10q+k+60}(s_1s_2^{-1}s_1s_2^{-1}s_1
  f^{5q-30})s_2^{-1}s_1s_2f^{10q+l-56} \\
&=s_2(s_1^{-1}s_2^{-1}s_1^{-1}s_2^{-1}s_1^{-1}s_2^{-1}s_1^{-1}
  s_2^{-1}s_1^{-1})s_2f^{5q+k+l-21} \\
&=s_2^{-1}s_1s_2s_1s_2^{-1}f^{4q+k+l-19}
\end{align*}
and go on by replacing $c''$
with $c'''=h''c''g''$, that is,
$$c'''=s_2^{-1}f^{4q+2k+2l-13}.$$
We further replace $g''$ and $h''$ by
\begin{align*}
g'''&=c'''g''{c'''}^{-1}=s_1f^{-10q+k+60} \\
\text{and}\qquad h'''&={c'''}^{-1}h''c'''=s_1f^{10q+l-56}.
\end{align*}
We lastly replace $h'''$ with $\bar h=g'''h'''=f^{k+l+5}$. These changes
clearly preserve the subgroup, thus $\langle g,h\rangle=\langle
g''',\bar h,c'''\rangle$. \fullref{intersectionwithf} now
implies that
$$\langle f\rangle\cap \langle g,h\rangle=\langle
\bar h, (g''')^2,(c''')^3,(g'''(c''')^{-1})^7\rangle$$
where
\begin{align*}
(g''')^2&=f^{-20q+2k+121}, \\
(c''')^3&=f^{12q+6k+6l-40} \\
\text{and}\qquad
(g'''(c''')^{-1})^7 &= f^{-97q-7k-14l+511}.
\end{align*}
This clearly implies the assertion of the lemma.
\end{proof}

Invertible generating pairs of this type however turn out to be extremely
rare, in fact there is only one manifold fibered over $S^2(2,3,7)$
that admits one, and for this manifold the invertible generating pair
is unique.

\begin{Lemma} Let $M$ be as above and
$$(g,h)=(s_1s_2^{-1}s_1s_2s_1s_2^{-1}s_1s_2^{-1}s_1s_2f^k,
s_2^{-1}s_1s_2s_1s_2^{-1}s_1s_2s_1s_2s_1f^l).$$
Then $(g,h)$ is an invertible generating pair if and only if $q=6$,
$k=-2$ and $l=-3$, that is,  $M$ is the Brieskorn homology sphere
$$\Sigma(2,3,7)= \{ z_{1}^{2} + z_{2}^{3} + z_{3}^{7} =0\} \cap \{ \vert
z_{1} \vert ^{2} + \vert  z_{2} \vert ^{2} + \vert z_{3} \vert ^{2}
=1\}.$$
Moreover this generating set is not geometric.
\end{Lemma}

\begin{proof}
Suppose first that $(g,h)$ is invertible. Let $\alpha$ be the
induced  automorphism of $G$. As the fibre is a power of $[g,h]$
and $\alpha([g,h])=[g^{-1},h^{-1}]$ is conjugate to $[g,h]$ it
follows that $\alpha$ preserves $f$. Thus $\alpha$ is induced by
a homeomorphism that is isotopic to the identiy, that is, $\alpha$
is an inner automorphism. Note that the inner automorphism must be
conjugation with $s_2s_1s_2^{-1}f^l$ for some $l$ by \fullref{inner237}
as $\alpha$ induces an involution on the base group. As the value of $l$
has no effect in the conjugation we can assume that $l=0$. Note that
$(s_2s_1s_2^{-1})^{-1}=s_2s_1^{-1}s_2^{-1}=f^{-1}s_2s_1s_2^{-1}$.
We obtain
\begin{align*}
f^{-1}s_2s_1s_2^{-1}\cdot g\cdot
s_2s_1s_2^{-1} &= f^{-1}s_2s_1s_2^{-1}\cdot s_1s_2^{-1}s_1s_2s_1
s_2^{-1}s_1s_2^{-1}s_1s_2f^k \cdot s_2s_1s_2^{-1} \\
&=s_2s_1s_2^{-1}s_1s_2^{-1}s_1s_2(s_1^{-1}s_2^{-1}s_1^{-1}
  s_2^{-1}s_1^{-1}s_2^{-1}s_1^{-1}s_2^{-1})f^{k+4} \\
&=s_2s_1s_2^{-1}s_1s_2^{-1}s_1s_2s_2s_1s_2s_1s_2s_1f^{-q+k+4} \\
&=s_2(s_1^{-1}s_2^{-1}s_1^{-1}s_2^{-1}s_1^{-1}s_2^{-1}s_1^{-1})
  s_2s_1s_2s_1f^{-q+k+9} \\
&=s_2(s_2s_1s_2s_1s_2s_1s_2)s_2s_1s_2s_1f^{-2q+k+9} \\
&=s_2^{-1}s_1^{-1}s_2s_1^{-1}s_2s_1^{-1}s_2^{-1}s_1^{-1}
  s_2s_1^{-1}f^{-2q+k+16}
\end{align*}
It follows that $f^{-1}s_2s_1s_2^{-1}\cdot g\cdot s_2s_1s_2^{-1}=g^{-1}$
if and only if $-k=-2q+k+16$, that is, if $k=q-8$. The analogous
calculation for $h$ shows that $f^{-1}s_2s_1s_2^{-1}\cdot h\cdot
s_2s_1s_2^{-1}=h^{-1}$ if and only if $-l=-2q+l-6$, that is, if $l=-q+3$.
 
\fullref{genset237} implies that $\langle g,h\rangle\cap
\langle f\rangle$ is generated by
\begin{align*}
f^{k+l+5} &= f^0, \\
f^{2k-20q+121} &=f^{2(q-8)-20q+121} \\ &=f^{-18q+105}, \\
f^{6k+6l+12q-40} &=f^{6(q-8)+6(-q+3)+12q-40} \\ &=f^{12q-70} \\
\text{and}\qquad
f^{-7k-14l-97q+511} &=f^{-7(q-8)-14(-q+3)-97q+511} \\ &=f^{-90q+525}.
\end{align*} Thus
$\langle g,h\rangle \cap\langle f\rangle=\langle f^{6q-35}\rangle$. Thus
$\langle g,h\rangle\cap\langle f\rangle=\langle f\rangle$ and therefore
$G=\langle g,h\rangle$ iff $q=6$ and therefore also $k=-2$ and $l=-3$.

To see that this generating set if not geometric it suffices to recall
that $M$ only admits two distinct Heegaard splitting, one vertical
and one horizontal corresponding to the generating pair discussed in
\fullref{case23r}. The generating pairs corresponding to these splittings
are not Nielsen equivalent to $(g,h)$ by \fullref{hyporb}. This implies
that $(g,h)$ is not geometric. The generating pair corresponding to the
horizontal splitting however happens  to induce the same involution of
the manifold.
\end{proof}

\section{Inversion on aspherical graph manifolds}\label{sec:involution}

In this section $M$ will be a closed, orientable, aspherical, sufficiently
complicated graph manifold, that is, it either admits  a non-empty splitting
along a finite collection of essential and pairwise non parallel tori into
Seifert fibered pieces or it is Seifert fibered with a hyperbolic base.
Such a decomposition is called a JSJ--decomposition. We show that any
inversion of $\pi_1(M)$ can be realized by an involution which preserves
the JSJ--splitting and the Seifert structures of the pieces. More precisely:

\begin{Proposition}\label{prop:involution} Let $M$ be a closed,
orientable, aspherical, sufficiently complicated graph manifold. Then
any inversion $\alpha$ of $\pi_1(M)$ can be realized by an involution
$\tau$ on $M$ with the following properties. \begin{enumerate}
\item The involution $\tau$ respects the JSJ--decomposition of $M$ and the Seifert fibered  structures on  the
JSJ--pieces. 
\item The underlying space of the quotient orbifold $M/\tau$ is
homeomorphic to $S^3$ and the map $p\co  M \to \vert M/\tau \vert = S^3$ is a
$2$--fold covering branched along a  prime and unsplittable link $L$ with
at most $3$--components. Moreover if $L$ has $\geq 2$ components either $L
= K_1 \cup K_2$ where $K_1$ is unknotted and $K_2$ is a 2--bridge knot, or
 $L= K_1 \cup K_2 \cup K_3$ where each component $K_i$ is unknotted and each sublink $K_i \cup K_j$,
$i \neq j \in \{1,2,3\}$ is a 2--bridge link.
\item The inversion $\alpha$ determines the pair $(S^3, L)$ up to
homeomorphism and is hyperelliptic if and only if $L$ is a 3--bridge link.
\end{enumerate}
\end{Proposition}

\begin{proof}
An inversion $\alpha \in \Aut(\pi_1(M))$ of the generating pair $(g,h)$
gives an extension of $\pi_1(M)$:
$$1 \to \pi_1(M) \to E \to \mathbb{Z}/2\mathbb{Z} \to 1,$$
such that $E= \langle \pi_1(M), t \, \, \vert  \, \, t^2,t \gamma t^{-1}
= \alpha (\gamma), \, \, \forall \gamma \in \pi_1(M) \rangle$.

The key step of the proof is to  realize $E$ as a group of orientation
preserving homeomorphisms of the universal cover $\wwtilde M$ of $M$
in such a way that the subgroup $\pi_1(M)$ is the deck transformation
group of the cover
$\wwtilde M \to M$. 

In the Seifert fibered case, using Kerckhoff's proof of the Nielsen's
conjecture \cite{Ker},  one can even realize $E$ as a group of fiber
preserving homeomorphisms of the universal cover $\wwtilde M$ of $M$,
according to Zimmerman and Zieschang~\cite{Zim1,ZZ}, see also \cite{Zie2}. 

In the case of a proper graph manifold this realization of $E$ as 
group of homeomorphisms of the universal cover $\wwtilde M$ is achieved
by Zimmerman~\cite[Section D]{Zim2}.

Now it is clear that the quotient group  $E/\pi_1(M) \cong
\mathbb{Z}/2\mathbb{Z}$ acts on $M$ and is generated by the projection
$\tau$ on $M$ of the  involution $t$ on $\wwtilde M$. Since $\wwtilde M$
is homeomorphic to $\mathbb{R}^3$, $t$ must have non-empty fixed point
set, hence $\tau$ has non-empty fixed point set.

The group $E$ is the orbifold fundamental group of the quotient
orbifold $\wwtilde M /E = M/\tau$. It is generated by the three
involutions $tg, th, t$, hence the fundamental group of the underlying
space $\vert M/\tau \vert$ is trivial (see Armstrong~\cite{Ar} and
Thurston~\cite[Chapter~13]{Thu}).  Since the orbifold $M/\tau$ is
either Seifert fibered or splits along a finite collection of disjoint
euclidean 2--suborbifolds into Seifert fibered pieces, it follows
that the underlying space $\vert M/\tau \vert$ must be homeomorphic
to the sphere $S^3$. Therefore $p\co  M \to \vert M/\tau \vert = S^3$
is a $2$--fold covering branched along a  link $L$ with at most $3$
components since $\dim H_1(M, \mathbb{Z}/2\mathbb{Z}) \leq 2$ because
$\pi_1(M)$ is 2--generated. Moreover $L$ is prime and unsplittable since
$M$ is irreducible.

The group $E$ is the so called $\pi$--orbifold group of the link $L$
(see Boileau and Zimmerman~\cite{BZim}). It can be computed from the
fundamental group $\pi_1(S^3\setminus L)$ of the link $L$ by killing the
squares of the meridian elements $\mu_i$:
$$ E = \pi_1(S^3\setminus L) /\langle\langle\mu_{i}^2\rangle\rangle.$$
The proof of the
Smith conjecture shows that a knot is unknotted if and only if its
$\pi$--orbifold group is $\mathbb{Z}/2\mathbb{Z}$. Moreover a prime and
unsplittable link is a  2--bridge link if and only if its $\pi$--orbifold
group is dihedral (see \cite{BZim}). The $\pi$--orbifold group of a
component $K$ of $L$ can be obtained from the $\pi$--orbifold group $E$
of $L$  by killing the meridians of the other components of $L$. Then the
assertion on the components of $L$ when  $L$ has at least $2$ components,
follows readily from the two previous facts.

By Zimmerman~\cite[Theorems~0.1 and 1.1]{Zim3}, the inversion $\alpha$
determines the involution $\tau$ up to conjugation by a homeomorphism
homotopic to the identity, hence $\alpha$ determines the link L. Therefore
$\alpha$ is hyperelliptic if and only if $\tau$ is induced by the
hyperelliptic involution of a Heegaard splitting of genus 2. By Birman
and Hilden~\cite{BH} this is equivalent to saying that the link $L$
has 3 bridges.
\end{proof}

The proof of \fullref{thm:seifert} follows from the results of
\fullref{sec:seifert} and \fullref{cor:seifert} below.

\begin{Corollary}\label{cor:seifert}
Let $M$ be an orientable Seifert manifold with base orbifold
$P^2(p,q)$ or $S^2(2,2,2,2l+1)$. Then any inversion of the fundamental
group of $M$ is hyperelliptic.
\end{Corollary}
\begin{proof}
Let $\alpha$ be an inversion of $\pi_1(M)$. By \fullref{prop:involution}
it is induced by a fiber preserving involution $\tau$. Let $p\co  M \to \vert
M/\tau \vert = S^3$ be the associated $2$--fold covering of $S^3$ branched along a  link $L$.

If the base space is  $S^2(2,2,2,2l+1)$, we apply the argument of Boileau
and Zieschang~\cite[Lemma 3.3]{BZ}. The link  $L$ has $3$ components since $\dim H_1(M, \mathbb{Z}/2\mathbb{Z}) = 2$. We distinguish two cases depending on the action of $\tau$ on the fibres:

\textbf{Case 1: $\tau$ reverses the orientation of the fibres.}\qua
Since the quotient space $\vert M/\tau \vert = S^3$ $\tau$ is a Montesinos
involution (see \cite{Mon1} and \cite[Chapter~4]{Mon2}). Thus $L$ is a
Montesinos link with $3$ components, one of which is a knotted 2--bridge
knot (see case (a) of \cite[Lemma 3.3]{BZ}). Hence $\tau$ cannot induce
an inversion.

\textbf{Case 2: $\tau$ preserves the orientation of the fibres.}\qua
Then Case 3.5 of the proof of \cite[Lemma 3.3]{BZ} shows that the quotient
space $\vert M/\tau \vert = S^3$ and $L$ has $3$ components, if and
only if $L$ is the 3--bridge torus link $T(3, 3(2l+1))$. Hence $\tau$
induces an inversion of $\pi_1(M)$ if and only if it is hyperelliptic.

When the base space is $P^2(p,q)$ then $\tau$ must reverse the orientation
of the singular fibres since the quotient space $\vert M/\tau \vert =
S^3$. Hence $\tau$ is a Montesinos involution  and $L$ is a generalized
Montesinos link (see \cite{Mon1} and \cite[Chapter~4]{Mon2}). It follows
from Boileau and Zieschang~\cite[Theorem~2.1]{BZ2} that the bridge number
of $L$ is $3$.  Therefore $\alpha$ is a hyperelliptic inversion.
\end{proof}

\begin{Remark}\label{rem:seifert}
The same type of arguments can be applied more generally to show that an
inversion of the fundamental group of a sufficiently complicated closed
orientable Seifert manifolds is hyperelliptic. The proof would follow
from the following two facts:
\begin{enumerate}
\item The bridge number  of a (generalized) Montesinos link $L$ with
$2$--fold branched cover $M$ is
$\rank\pi_1(M) + 1$, except in the case when $M$ has base space $S^2(2,2,\ldots,
2,2l+1)$ with an even number of singular points, see \cite{BZ2}.
\item The bridge number of a torus link whose $\pi$--orbifold group is
generated by $k$ elements of order $2$ is at most $k$.
\end{enumerate}

However the interest of the algebraic proof given in the previous \fullref{sec:seifert} is to give more precise results: it shows that many Heegaard genus 2 Seifert manifolds have non invertible generating pairs and even that the Brieskorn 
homology sphere admits an invertible generating pair which is not geometric.
\end{Remark}

\section{Graph manifolds}\label{sec:graph}

\fullref{thm:seifert} clearly implies \fullref{main} for Seifert manifolds. In this section we consider 
graph-manifolds which are not Seifert fibered, that is, graph-manifolds that
have a non-trivial JSJ--decomposition.
 
 The purpose of the section is to establish the following proposition, which finishes the proof of  
\fullref{main} together with \fullref{thm:seifert}.

\begin{Proposition}\label{prop:genus3} The fundamental group of a closed, orientable graph manifold of 
Heegaard genus greater than $ 2$ does not admit an inversion.
\end{Proposition}
 
We need only to consider the case where the graph manifold is not Seifert fibered.
The following lemma precisely describes such graph manifolds of Heegaard genus $3$. 

\begin{Lemma}[Weidmann
\cite{We}]\label{lem:cases} Let $M$ be a graph manifold of Heegaard genus
$3$ which is not Seifert fibered and has $2$--generated fundamental group. Then the following hold:\begin{enumerate}

\item $M=M_1\cup _{T}M_2$ where $M_1$ is fibered over M\"o or M\"o$(p)$ and $M_2$ is fibered over $D(2,2l+1)$. 
\item The intersection number of the the fibres $f_{1}$ and $f_{2}$ on $T$  is 
$\Delta(f_{1}, f_{2}) = \pm 1$.
\item The pair  $(M_1,f_{2})$ is not the exterior of a 1--bridge knot in a lens space with  $f_{2}$ being a meridian.
\item The pair  $(M_2,f_{1})$ is not the exterior of a 2--bridge knot in $S^3$ with  $f_{1}$ being a meridian.
\end{enumerate}
\end{Lemma}

\begin{proof}[Proof of \fullref{prop:genus3}] 
If $\pi_1(M)$ admits an inversion  it follows from
\fullref{prop:involution} that it is induced by an involution $\tau$,
unique up to conjugacy and which preserves each Seifert piece and their
Seifert fibrations.  We distinguish the cases where the involution $\tau$
has fixed points on the JSJ--splitting torus  $T$ and where it does not.

\textbf{Case 1: $\Fix(\tau) \cap T \not = \emptyset$}\qua Then
$\tau$ induces an hyperelliptic involution on $T$ and the fixed point set
of $\tau$ meets the splitting  torus $T$ in $4$ points. The quotient
spaces $M_{1}/\tau$ and $M_2/\tau$ together with the image of the fixed
point set of $\tau$ are two non-rational Montesinos tangles (see Bonahon
and Siebenmann~\cite{BS}): we use notations analogous to the ones used for
Montesinos links in Boileau and Zieschang~\cite[Section 2]{BZ2} or Burde
and Zieschang~\cite[Chapter~12]{BuZ}

$B_1 = M_{1}/\tau = B(-1; \beta/p)$ where $\beta/p$ is the type of the
singular fiber of order $p$ of $M_{1}$, see \fullref{fig1}.

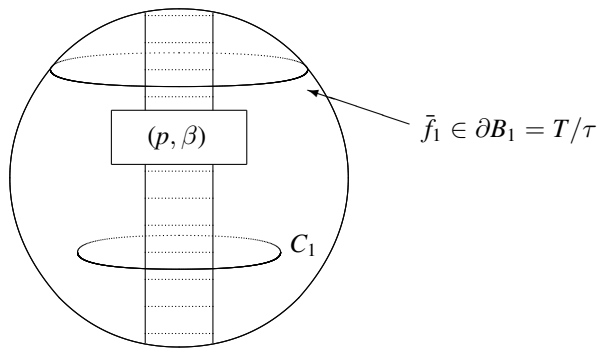
\begin{figure}[ht!]
  \centerline{ \footnotesize \setlength{\unitlength}{.9cm}
\begin{picture}(15,5)
\small
\bezier{300}(7.5,0)(8.55,0)(9.27,.73)
\bezier{300}(7.5,0)(6.45,0)(5.73,.73)
\bezier{300}(7.5,5)(8.55,5)(9.27,4.27)
\bezier{300}(7.5,5)(6.45,5)(5.73,4.27)
\bezier{300}(10,2.5)(10,1.55)(9.27,.73)
\bezier{300}(5,2.5)(5,1.55)(5.73,.73)
\bezier{300}(10,2.5)(10,3.45)(9.27,4.27)
\bezier{300}(5,2.5)(5,3.45)(5.73,4.27)
\put(6.5,2.7){\line(1,0){2}}\put(7,3){$(p, \beta)$}
\put(6.5,3.5){\line(1,0){2}}
\put(6.5,3.5){\line(0,-1){.8}}
\put(8.5,3.5){\line(0,-1){.8}}
\put(8,3.5){\line(0,1){1.44}}
\put(7,3.5){\line(0,1){1.44}}
\put(8,2.7){\line(0,-1){2.64}}
\put(7,2.7){\line(0,-1){2.64}}
\bezier{50}(5.6,4.1)(5.6,4.35)(7.5,4.35)
\bezier{50}(9.4,4.1)(9.4,4.35)(7.5,4.35)
\bezier{500}(5.6,4.1)(5.6,3.85)(7.5,3.85)
\bezier{500}(9.4,4.1)(9.4,3.85)(7.5,3.85)
\put(10.9,3.3){\vector(-3,1){1.5}}\put(11.1,3.1){$\bar f_1\in\partial B_1=T/\tau$}
\bezier{50}(6,1.4)(6,1.65)(7.5,1.65)
\bezier{50}(9,1.4)(9,1.65)(7.5,1.65)
\bezier{500}(6,1.4)(6,1.15)(7.5,1.15)
\bezier{500}(9,1.4)(9,1.15)(7.5,1.15)
\put(9.1,1.4){$C_1$}
\bezier{20}(7,.6)(7.5,.6)(8,.6)
\bezier{20}(7,.2)(7.5,.2)(8,.2)
\bezier{20}(7,1)(7.5,1)(8,1)
\bezier{20}(7,1.8)(7.5,1.8)(8,1.8)
\bezier{20}(7,1.4)(7.5,1.4)(8,1.4)
\bezier{20}(7,2.2)(7.5,2.2)(8,2.2)
\bezier{20}(7,2.6)(7.5,2.6)(8,2.6)
\bezier{20}(7,3.7)(7.5,3.7)(8,3.7)
\bezier{20}(7,4.1)(7.5,4.1)(8,4.1)
\bezier{20}(7,4.5)(7.5,4.5)(8,4.5)
\bezier{20}(7,4.9)(7.5,4.9)(8,4.9)
\end{picture}}
\caption{$B_1=B(-1\,|\,\beta/p)$}
\label{fig1}
\end{figure}
 
$B_2 = M_2/\tau = B(0; 1/2, \lambda/2l+1)$ where $\lambda/2l+1$ is the
type of the singular fiber of order $2l+1$ of $M_{2}$, see \fullref{fig2}.

\begin{figure}[ht!]
  \centerline{ \footnotesize \setlength{\unitlength}{.8cm}
\begin{picture}(15,5)
\small
\bezier{300}(7.5,0)(8.55,0)(9.27,.73)
\bezier{300}(7.5,0)(6.45,0)(5.73,.73)
\bezier{300}(7.5,5)(8.55,5)(9.27,4.27)
\bezier{300}(7.5,5)(6.45,5)(5.73,4.27)
\bezier{300}(10,2.5)(10,1.55)(9.27,.73)
\bezier{300}(5,2.5)(5,1.55)(5.73,.73)
\bezier{300}(10,2.5)(10,3.45)(9.27,4.27)
\bezier{300}(5,2.5)(5,3.45)(5.73,4.27)
\put(6.5,2.5){\line(1,0){2}}
\put(6.5,2.5){\line(0,-1){.8}}
\put(8.5,2.5){\line(0,-1){.8}}
\put(8,1.7){\line(0,-1){1.64}}
\put(7,1.7){\line(0,-1){1.64}}
\put(8,4.94){\line(0,-1){.84}}
\put(7,4.94){\line(0,-1){.84}}
\put(8,2.5){\line(0,1){.7}}
\bezier{300}(8,4.1)(8,3.5)(7.5,3.5)
\bezier{300}(7.25,3.6)(7.35,3.5)(7.5,3.5)
\bezier{300}(7,4.1)(7,3.95)(7.1,3.75)
\bezier{300}(8,3.2)(8,3.45)(7.9,3.55)
\bezier{300}(7,3.2)(7,3.8)(7.5,3.8)
\bezier{300}(7.75,3.7)(7.8,3.7)(7.5,3.8)
\put(7,2.5){\line(0,1){.7}}
\bezier{20}(7,.6)(7.5,.6)(8,.6)
\bezier{20}(7,.2)(7.5,.2)(8,.2)
\bezier{20}(7,1)(7.5,1)(8,1)
\bezier{20}(7,1.4)(7.5,1.4)(8,1.4)
\bezier{20}(7,2.8)(7.5,2.8)(8,2.8)
\bezier{20}(7,3.2)(7.5,3.2)(8,3.2)
\bezier{20}(7,4)(7.5,4)(8,4)
\bezier{20}(7,4.4)(7.5,4.4)(8,4.4)
\bezier{20}(7,4.8)(7.5,4.8)(8,4.8)
\bezier{6}(7.5,3.5)(7.5,3.65)(7.5,3.8)
\bezier{300}(5.3,3.65)(5.3,3.25)(7.5,3.25)
\bezier{300}(9.7,3.65)(9.7,3.25)(7.5,3.25)
\bezier{50}(5.3,3.65)(5.3,4.05)(7.5,4.05)
\bezier{50}(9.7,3.65)(9.7,4.05)(7.5,4.05)
\put(11,2.7){\vector(-3,1){1.5}}\put(11.1,2.6){$\bar f_2\in\partial B_2=T/\tau$}
\scriptsize\put(6.5,1.7){\line(1,0){2}}\put(6.7,2){$(2l+1, \lambda)$}
\end{picture}}
\caption{$B_2=B(0\,|\,1/2,\lambda/(2l+1))$}
\label{fig2}
\end{figure}
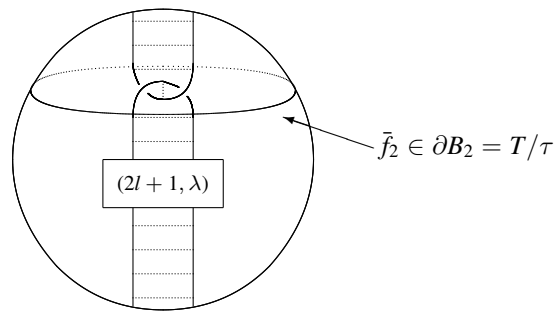

The manifold $M$ is a $2$--fold covering of $S^3$ branched along a so called Conway algebraic link $L$ obtained by gluing 
the two tangles $B_1$ and $B_2$ along their boundary (see Bonahon and
Siebenmann~\cite{BS}).
The link  $L =K \cup C_1$, where $K$ is a 2--bridge knot or  link by
\fullref{prop:involution}, and $C_1$ is an unknotted component of the
Montesinos tangle $B_1$ which corresponds to a $S^1$--fiber of the orbifold
Seifert  fibration of $B_1$ obtained by quotienting the Seifert fibration
of $M_1$ (see \fullref{fig1} and 
Montesinos~\cite{Mon1,Mon2}).
We remark that when forgetting the component $C_1$, the Montesinos tangle $B_1$ becomes a rational tangle 
$B'_1$ of type $\beta/p \in \mathbb{Q}/\mathbb{Z}$ with the orbifold Seifert fibration induced from the one of $B_1$.
The rational tangle $B'_1$ contains a properly embedded meridian disk 
$(\Delta_1, \partial \Delta_1) \hookrightarrow (B'_1, \partial B'_1)$ which separates its two strands. Let 
$\delta_1 = \partial \Delta_1$ be the simple closed curve on the pillow  $\partial B'_1 = \partial B_1$. 
In order to obtain a $2$--bridge knot or link 
$K$ when gluing the rational tangle $B'_1$ to the Montesinos tangle $B_2$,  one needs that the 
orbifold Seifert fibration of $B_2$ (inherited from  the Seifert fibration
of $M_2$) extends to an orbifold Seifert fibration on the tangle $B'_1$
without singular fiber. Otherwise the $2$--fold branched covering of $K$ obtained by Dehn filling $M_2$ would have three singular fibers and could not be a Lens space.

Let $\bar f_1$ and $\bar f_2$ be the images on 
$\partial B'_1 = \partial B_2$ of the fibers $f_1 \subset \partial M_1 = T$ and $f_2 \subset \partial M_2 = T$.
We say that two simple closed curves on the pillow $\partial B'_1 = \partial B_2$ are dual if they are the images of two 
simple closed curves on the torus $T$ which are dual, that is, meets in one point. Since $\Delta(f_{1}, f_{2}) = \pm 1$ it follows that $\bar f_1$ and $\bar f_2$ are dual on the pillow $\partial B'_1 = \partial B_2$, and thus the orbifold Seifert fibration of the pillow $\partial B'_1 = \partial B_2$ with fiber $\bar f_2$ extends to an orbifold Seifert  fibration of the tangle $B'_1$ with a singular fiber of order $\beta$. Hence $\beta = \pm 1$. In this case $\bar f_2$ is also dual to the 
curve $\delta_1 \subset \partial B'_1 = \partial B_2$. Therefore, up to orientation, there is at most one possibility for choosing $\bar f_2$ on the pillow $\partial B'_1 = \partial B_2$, and thus for choosing $f_2$ on $T$.

The condition $\beta = \pm 1$ is equivalent for $M_1$ to be the exterior
of a one-bridge knot in a Lens space by Weidmann~\cite[Lemma 1]{We}. Therefore the existence of an inversion on $\pi_1(M)$ implies that $M_1$ is the exterior of a one-bridge knot in a lens space and that, up to orientation of $M_1$, there is at most one possible choice for the fiber $f_2$ on 
$T= \partial M_1$. Kobayashi's construction of Heegaard genus two graph-manifolds shows that this happens 
precisely when $f_2$ is the meridian of the knot space in the Lens space, which is the excluded case.

\textbf{Case 2: $\Fix(\tau) \cap T = \emptyset$}\qua In this case
the quotient spaces $V_1= M_{1}/\tau$ and $V_2=M_2/\tau$ are Seifert
fibered 3--orbifolds with boundary the torus $T/\tau$ such that $V_1 \cup V_2 = S^3$. 

The action induced by $\tau$ on the base $D(2,2l+1)$ of $M_2$ must be the identity, hence the restriction of
$\tau$ to ${M_{2}}$ is the $\mathbb{Z}/2\mathbb{Z}$ action embedded into
the $S^1$--action defining the Seifert fibration on $M_2$. In particular
the exceptional fibre of order $2$ of $M_2$ is fixed pointwise by $\tau$.
Thus $V_2$ is a Seifert fibered solid torus with base $D(2l+1)$ and the
component $K_2 \subset V_2$ of the branching link L  is a regular fiber in
this Seifert fibered solid torus, by Seifert~\cite[Section~14]{Sei}. 

Since the quotient $M_1/\tau$ embeds in $S^3$ with boundary a torus, the action induced by $\tau$ on the base M\"o or M\"o($p$) cannot be the identity and cannot have fixed point on the boundary, hence it is the reflection through the core of the underlying M\"obius strip, where the singular point of order $p$ is assumed to lie on this core.

 To understand the restriction of the involution $\tau$ to $M_1$ we first
consider the case where $M_1$ is the twisted $I$--bundle over the Klein
bottle, that is, the base is a M\"obius strip. In this case the quotient $V_1=
M_{1}/\tau$  is the Euclidean 3--orbifold which is an orientable bundle (in
the orbifold sense) over a M\"obius strip with silvered boundary, where a
boundary is called silvered if it corresponds to a reflection axis (see,
for example, Boileau and Porti~\cite[Chapter~4]{BP}). The underlying space is a solid 
torus which is a neighborhood of the M\"obius strip and the branching
locus $K_1$ is the boundary of the M\"obius strip. The quotient $V_1$ is a
Seifert 3--orbifold: each meridian disk of the solid torus $V_1$ is fibered by concentric circles together with an arc joining the two intersection points of the branching locus with the meridian disk.

To understand the restriction of the involution $\tau$ to $M_1$ when the
base is M\"o$(p)$, we have just to pull out a tubular neighborhood of a
$\tau$--invariant regular fiber in the previous example and glue it back to
insert a singular fiber of type $\beta/p$. In the quotient $V_1$ it
corresponds to pulling out a tubular neighborhood of a fiber which is an
interval with silvered boundary (that is, a non-separating arc of the M\"obius 
strip bounded by the branching locus $K_1$) and inserting a rational tangle of type $\beta/p$. The quotient 
$V_1$ is still a Seifert 3--orbifold with underlying space a solid torus
$V_1$ and whose boundary $\bar T = T/\tau$ is fibered by the boundaries of
meridian disks of $V_1$, see \fullref{fig3}.

\begin{figure}[ht!]
  \centerline{ \footnotesize \setlength{\unitlength}{.8cm}
\begin{picture}(15,6)
\small
\thicklines
\bezier{300}(7.5,0)(9.6,0)(11.04,.73)
\bezier{300}(7.5,0)(5.4,0)(3.96,.73)
\bezier{300}(7.5,5)(9.6,5)(11.04,4.27)
\bezier{300}(7.5,5)(5.4,5)(3.96,4.27)
\bezier{300}(12.5,2.5)(12.5,1.55)(11.04,.73)
\bezier{300}(2.5,2.5)(2.5,1.55)(3.96,.73)
\bezier{300}(12.5,2.5)(12.5,3.45)(11.04,4.27)
\bezier{300}(2.5,2.5)(2.5,3.45)(3.96,4.27)
\bezier{300}(7.5,2.3)(8.8,2.3)(9.5,2.7)
\bezier{300}(7.5,2.3)(6.2,2.3)(5.5,2.7)
\bezier{300}(7.5,2.8)(6.5,2.8)(6,2.53)
\bezier{300}(7.5,2.8)(8.5,2.8)(9,2.53)\thinlines
\bezier{300}(8.1,.8)(6.5,.80)(5.5,1.12)
\bezier{300}(3.6,2.5)(3.6,1.8)(5.5,1.12)
\bezier{300}(8.1,4.2)(6.5,4.2)(5.5,3.88)
\bezier{300}(3.6,2.5)(3.6,3.2)(5.5,3.88)
\bezier{300}(8.1,1.5)(6.8,1.5)(6,1.72)
\bezier{300}(4.6,2.5)(4.6,2.1)(6,1.72)
\bezier{300}(8.1,3.5)(6.8,3.5)(6,3.28)
\bezier{300}(4.6,2.5)(4.6,2.9)(6,3.28)
\bezier{300}(8.1,3.5)(9,3.5)(9.5,3.38)
\bezier{300}(11.4,2.5)(11.4,3)(9.5,3.38)
\bezier{300}(11.4,2.5)(11.4,1.3)(8.9,.8)
\bezier{300}(8.1,4.2)(9.6,4.2)(10.2,3.4)
\bezier{300}(8.9,1.5)(9.5,1.5)(10,1.7)
\bezier{300}(10,1.7)(11.2,2.3)(10.5,3)
\bezier{300}(7.2,3.47)(7.15,3.8)(7.1,4.16)
\bezier{300}(6.7,.04)(7,.04)(7.1,1.1)
\bezier{300}(6.9,2.3)(7.2,2.29)(7.1,1.1)
\bezier{200}(6.7,.04)(6.5,.04)(6.45,.88)
\bezier{200}(6.9,2.3)(6.6,2.31)(6.5,1.6)
\bezier{10}(6.45,.88)(6.475,1.24)(6.5,1.6)
\put(8.5,1.15){\circle{1.1}}\put(8.1,1.08){$\beta/p$}
\put(11.3,4.6){\vector(-3,-1){1.5}}\put(11,4.8){$\hbox{Branching locus }K_1$}
\put(5.7,5.2){\vector(1,-1){1.3}}\put(3,5.4){$\hbox{Fibre with silvered boundary}$}
\put(4.4,.1){\vector(4,1){2}}\put(2.5,0){$\bar f_1\subset \partial V_1$}
\end{picture}}
\caption{The Seifert fibered orbifold $V_1$}
\label{fig3}
\end{figure}
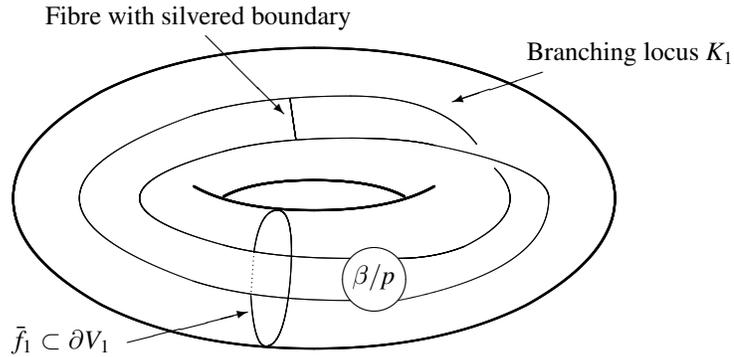

Let $L= K_1 \cup K_2$ be the branching link $L \subset S^3$, with  $K_1 \subset V_1$ and $K_2 \subset V_2$. 
Let $\Delta_1$ and $\Delta_2$ be  meridian disks of the solid tori $V_1$ and $V_2$ respectively. Since  $V_1 \cup V_2$ is $S^3$, their boundaries verify:  
$\Delta (\partial \Delta_1, \partial \Delta_2) = \pm 1$. Thus $\partial \Delta_1$ and $\partial \Delta_2$ form a base of $H^1(\bar T,  \mathbb{Z})$ for the torus quotient $\bar T = T/\tau$. Let $\bar f_1$ and $\bar f_2$ be the images on $\bar T$ of the fibres $f_1$ and $f_2$ on $T$. 

Since $K_2$ is isotopic to the regular fibre $\bar f_2$ of the fibered solid torus $V_2$, it can be pushed onto
$\bar T = \partial V_1 = \partial V_2$. Moreover $K_2 = \bar f_2$ will meet $f_1 = \partial \Delta_1$ in one point, since $f_2$ meets $f_1$ geometrically in one point on $T= \partial M_1 = \partial M_2$. Therefore 
$K_2$ is unknotted in $S^3$ and the 2--fold covering of $V_2$ branched over
$K_2$ embeds into the 2--fold branched covering of $S^3$ branched along
$K_2$, which is still $S^3$. Hence $M_2$ is the exterior of a torus knot
of type $(2,2l+1)$, which is a $2$--bridge knot.

Let $\mu_2$ be the meridian curve on $\partial M_2$ when $M_2$ is considered as a knot exterior in $S^3$. 
The Seifert surface $F_2$ of the fibred torus knot exterior $M_2$ verifies that  $\Delta(\mu_2, \partial F_2) = 1$ and $\Delta(\partial F_2, f_2) = \pm 2(2l+1)$. Moreover $\partial F_2$ projects to 
the boundary $\partial \Delta_2$ of the meridian disk of $V_2$.  

As the fiber $\bar f_1 = \partial \Delta_1$ meets $\partial \Delta_2$ in one point, we have $\vert \Delta (f_1, \partial F_2) \vert \leq 2$. Therefore 
we have shown that:

\begin{enumerate}
\item $\vert \Delta(\partial F_2, f_2)\vert =  2(2l+1)$
\item $\vert \Delta (f_1, \partial F_2) \vert \leq 2$
\item $\vert \Delta(f_1, f_2) \vert = 1$
\end{enumerate}

By (1), $f_2 = \partial F_2 \pm 2(2l+1) \mu_2$ on $\partial M_2 = T$.

By (2) and (3)  $2(2l+1) \vert \Delta(\mu_2, f_1)  \vert \leq 3$, and hence $\Delta(\mu_2, f_1) = 0$.
This would imply that $\mu_2$ and $f_1$ are isotopic on $T = \partial M_2$ contradicting our hypothesis.

Another argument to get a contradiction would be  to show that, after
pushing $K_2$ into the unknotted solid torus $V_1$, $L$ has at most 3
bridges using the method in Boileau and Zieschang~\cite{BZ2}.
\end{proof}

\bibliographystyle{gtart}
\bibliography{link}

\end{document}